\documentclass[12pt]{article}
\usepackage[amsmath]{e-jc}

\usepackage{amsmath,amsthm,amssymb}
\usepackage{enumitem}
\usepackage{graphicx}
\usepackage{tikz}
\usepackage{diagbox}
\ifpdf%
\usepackage{epstopdf}%
\else%
\fi
\newcommand{\Ftwobar}{\overline{\mathbb{F}_2}}
\specs{n}{vol (issue)}{year}{}
\dateline{Nov 8, 2023}{TBD}{TBD}

\title{Ramsey numbers through the lenses of polynomial ideals and Nullstellens\"atze}

\author{Jes\'us A. De Loera\authornote{1}
\and
William J. Wesley\authornote{2}
}

\authortext{1}{Department of Mathematics, University of California, Davis
   (\email{deloera@math.ucdavis.edu}).}
\authortext{2}{Department of Mathematics, University of California, San Diego (\email{wjwesley@ucsd.edu}).}

\Copyright{The authors. Released under the CC BY-ND license (International 4.0).}

\begin{document}

\maketitle
\begin{abstract}
We study Ramsey numbers via Hilbert's Nullstellensatz and Alon's Combinatorial Nullstellensatz. We give encodings of various types of Ramsey 
numbers, from the classical graph theory version to the arithmetic Ramsey numbers of Rado, Schur, and van der Waerden, as systems of polynomial equations whose solutions are in bijection to colorings that avoid monochromatic patterns. For example, in the classical graph theory case the solutions correspond to Ramsey graphs of order $n$, those that do not contain a copy of $K_r$ or $\overline{K_s}$. 
When these systems of equations have no solution for the first time, the Ramsey-type number in question is attained. We construct Hilbert Nullstellensatz certificates whose degrees are equal to the restricted online Ramsey numbers introduced by Conlon, Fox, Grinshpun and He. Similar results apply to many numbers in Ramsey theory, including Rado, van der Waerden, Schur, and Hales-Jewett numbers. Finally, inspired by work of Alon and Tarsi, we introduce a new family of numbers that relate to the coefficients of a certain ``Ramsey polynomial" and give lower bounds for Ramsey numbers. Our work reveals connections to the computational complexity of Ramsey numbers.
\end{abstract}

\newcommand\blfootnote[1]{%
  \begingroup
  \renewcommand\thefootnote{}\footnote{#1}%
  \addtocounter{footnote}{-1}%
  \endgroup
}
\blfootnote{A less developed version of this paper appeared in \cite{WJW_FPSAC}.}
\section{Introduction} 
Ramsey numbers are some of the most interesting and mysterious combinatorial numbers \cite{GrahamRothschildSpencer}. They appear in graph theory, geometry, number theory, 
and other fields. In this paper we use Hilbert's Nullstellensatz to uncover a new algebraic dependence between Ramsey numbers and restricted online Ramsey numbers \cite{RestrictedOnlineRamseyNumbers_Conlon_et_al}. This new dependence applies for all types of Ramsey-numbers (arithmetic, geometric, graph theoretic). We also use a new ``Ramsey polynomial", in the spirit of Alon's combinatorial Nullstellensatz \cite{AlonCombinatorialNullstellensatz}, to show lower bounds on Ramsey numbers come from the vanishing of its coefficients. Overall, we point to algebraic ways 
to measure the computational complexity of Ramsey numbers.

In its simplest, most popular form, the graph theory \textit{Ramsey number} $R(r,s)$ is the smallest positive integer $n$ such that every 2-coloring of the edges of $K_n$ contains a 
copy of $K_r$ in the first color or $K_s$ in the second color. Ramsey numbers can be generalized by allowing more than two colors and graphs other than $K_r$ and $K_s$. The number $R(G_1,G_2,\dots,G_k)$ is the smallest positive integer $n$ such that every $k$-coloring of the edges of $K_n$ contains a copy of $G_i$ in color $i$. If $G_i = K_{r_i}$ for all $i$, we simply write $R(r_1,r_2,\dots,r_k)$. All of these numbers are finite by Ramsey's theorem \cite{Ramsey}.  


In the concrete special case of graph theory, we introduce two interpretations of $R(r,s)$ and related numbers and graphs using polynomial ideals, varieties, and Hilbert Nullstellensatz identities (see \cite{CoxLittleOShea} for an introduction). 
We tie the values and complexity of $R(r,s)$ to encodings or models using systems of  \emph{polynomial equations}. We show that these encodings give interesting information on the computational complexity of Ramsey numbers and Ramsey graphs. The precise models appear in Section \ref{PolynomialSection}. Before we state our main results, let us recall some relevant prior context and results:

Computing exact values for Ramsey numbers is a challenge. In fact, there are only 
nine values of $R(r,s)$ with $3 \le r \le s$ whose exact values are currently known, and the only known non-trivial Ramsey numbers with more than two colors are $R(3,3,3) = 17$ 
and $R(3,3,4)$ = 30 \cite{R333Equals17,R334Equals30}. Ramsey numbers as small as $R(5,5)$ remain unknown, and the best known bounds are $43\le R(5,5) \le 46$. See \cite{ExooR55LowerBound,R55Le46}. The numbers $R(G_1,G_2)$ 
are known for some families of graphs, but many cases remain open (see, for example, \cite{RamseySurvey} for a survey of small Ramsey numbers and their best-known bounds). The best known asymptotic lower bounds for diagonal Ramsey numbers $R(s,s)$ are given in \cite{SpencerRamseyLB}, and the best upper bounds are from \cite{gupta2024optimizingcgmsupperbound}, following the breakthrough exponential improvement in \cite{campos2025exponentialimprovementdiagonalramsey}. 

Although we know that in practice computing Ramsey numbers is extremely difficult (and considered harder than fighting a war with an alien civilization), it is not clear what is the appropriate computational complexity class to show theoretical hardness of computing Ramsey numbers $R(r,s)$. For example, the closely related \emph{arrowing decision problem} asks whether, given three graphs $F, G, H$ is there is a red/blue edge-coloring of $F$ that contains either a red $G$
or a blue $H$? This decision problem was shown to be in co-NP for fixed choices of $G,H$ \cite{BurrRamseyComplexity}. Later Schaefer \cite{SchaeferRamseyComplexity} showed that in general it is in the polynomial hierarchy to answer these queries, but it is not clear what to do with this complexity question when $F,G,H$ are complete graphs $K_N,K_r,K_s$ because there is only one value $R(r,s)$ for each input $N,r,s$, thus it is not clear how it can be hard for any of the usual classes like NP. See details in \cite{SchaeferRamseyComplexity,BurrRamseyComplexity, Haanpaa00computationalmethods}. 

In recent years, Pak and collaborators \cite{PakComplexityEnumCombinatorics,Ikenmeyer+Pak,PakCombinatorialInterpretation} have proposed that another way to measure computational complexity is by looking at counting sequences. Here we propose that their point of view could be another way to assert hardness of $R(r,s)$ by counting of Ramsey graphs: \emph{Ramsey $(r,s)$-graphs}, are graphs with no clique of size $r$, and no independent set of size $s$. Clearly, the number of vertices of a Ramsey $(r,s)$-graph is less than the Ramsey number $R(r,s)$. We are interested in the number of Ramsey graphs on $n$ vertices denoted by $RG(n,r,s)$. What is the complexity of counting the sequence of numbers $\{RG(n,r,s)\}_{n=1}^\infty$? From Ramsey's theorem this sequence consists of $R(r,s)-1$ positive numbers and then an infinite tail of zeroes. We give some examples of $RG(n,r,s)$ in Table \ref{TableRG}, which are computed using the \#SAT solver {\scshape relsat} \cite{Relsat}, following the work of successful computations in Ramsey theory using SAT solvers in, for example, \cite{SchurFive} and \cite{WJW_Rado_ISSAC} (see also the database at  \cite{McKayRamseyGraphs}). Thus, the hardness of $R(r,s)$ can then be rephrased as the following question:
\vskip 12pt
{\bf QUESTION:} Is the counting function  $RG(n,r,s)$ in $\#P$? 

\begin{table}[h]
\caption{Table of small values of function $RG(n,r,s)$ counting Ramsey $(r,s)$-graphs.}
\label{TableRG}
\begin{center}
	\begin{tabular}{c|c|c|c|c|c|c}
		$n$&$RG(n,3,3)$&$RG(n,3,4)$&$RG(n,3,5)$ & $RG(n,3,6)$ & $RG(n,3,7)$ & $RG(n,4,4)$\\
		\hline
		1 & 1 & 1 & 1 & 1 & 1 & 1\\
		2 & 2 & 2 & 2 & 2 & 2 & 2\\
		3 & 6 & 7 & 7 & 7 & 7 & 8 \\
		4 & 18 & 40 & 41 & 41 & 41 & 62\\
		5 & 12 & 322 & 387 & 388 & 388 & 892\\
		6 & 0 & 2812 & 5617 & 5788 &5789 & 22484\\
		7 & 0 & 13842 & 113949 & 133080 & 133500 & 923012\\
		8 & 0 & 17640 & 2728617 & 4569085 & 4681281 & 55881692
		\\
		9 & 0 & 0 & 55650276 & 220280031 & 245743539 & 4319387624\\
	\end{tabular}
\end{center}
\end{table}


Most of the results of this paper are indeed motivated by understanding and counting the Ramsey graphs as the solutions of a system of polynomial equations. 

\subsection*{Our Contributions:}

Our first contribution, in Section  \ref{PolynomialSection}, is to reintroduce the sequence $\{RG(n,r,s)\}_{n=1}^\infty$  as the number of solutions of certain zero-dimensional ideals over
the polynomial ring $\overline{\mathbb{F}_2}[x_1,\dots,x_n].$ The solutions are indicator vectors that yield all Ramsey graphs (note, here they are not counted up to symmetry or automorphism classes). Some simple properties of $RG(n,r,s)$, such as the fact that $RG(n,r,s) \le RG(n,r+1,s)$, follow immediately from Theorem \ref{binEdge2}.

\begin{theorem}\label{binEdge2}
	The Ramsey number $R(G_1,\dots,G_k)$ is at most $n$ if and only if there is no solution to the following system over $\overline{\mathbb{F}_2},$ where $K_n=(V,E)$ is the complete graph on $n$ vertices. Moreover, when the system has solutions, the number of solutions to this system is equal to the number of graphs of order $n$ that avoid copies of $G_i$ in color $i$. In particular, when $k=2$, $G_1 = K_r$, and $G_2 = K_s$, this is the number of Ramsey graphs $RG(n,r,s)$.
	
	\begin{align} \label{thebest}
	p_{H,i} := \prod_{e\in E(H)} x_{i,e} &=0 &&{\forall i, 1\le i \le k,\quad \forall  H\subseteq K_n,\ H \cong G_i},  \\
	q_e := 1+\sum_{i=1}^k x_{i,e}&= 0 &&{\forall e \in E}, \\
	u_{i,j,e} := x_{i,e}x_{j,e} & = 0 &&\forall e\in E,\ \forall i,j,\ i\neq j.
	\end{align}
	
	When $k = 2$, $G_1 = K_r$, and $G_2 = K_s$, the \emph{Ramsey ideal} $RI(n,r,s)$ is ideal of the polynomial ring $\Ftwobar[ x_{1,e},x_{2,e}]_{e\in E(K_n)}$ generated by the polynomials $p_H$, $q_e$ and $u_{i,j,e}$. Then we have 
	
	
	
	$$RI(n,r,s) \supseteq RI(n,r+1,s) \supseteq  \dots \supseteq RI(n,n,s) \supseteq RI(n,n+1,s) = RI(n,n+2,s) = \dots$$
	and 
	$$RI(n,r,s) \supseteq RI(n,r,s+1) \supseteq  \dots \supseteq RI(n,r,n) \supseteq RI(n,r,n+1) = RI(n,r,n+2) = \dots$$
	
	The very first value of $n$ for which the system of equations has no solution is equal to the Ramsey number.
\end{theorem}






One important result for our purposes is the famous \emph{Hilbert's Nullstellensatz}, which states that a system of polynomial equations $f_1 =\dots = f_m = 0$ over an algebraically closed field $K$ has no solution if and only if there exist coefficient polynomials $\beta_1,\dots,\beta_m$ such that \begin{equation}\label{NullCert}
\sum_{i=1}^m \beta_i f_i = 1.
\end{equation}

We call such an identity \eqref{NullCert} a \textit{Nullstellensatz certificate}. The \textit{degree} of a certificate is the largest degree of the polynomials $\beta_i$. Note that in our case the existence of a Nullstellensatz certificate is equivalent to an upper bound on the Ramsey number. The strong connection between combinatorial problems and the Hilbert Nullstellensatz has been investigated in, for example,  \cite{BussPitassi_InductionPrinciple,MarguliesOnnPasechnik_Nullstellensatz_Partition,deloera_lee_margulies_onn_2009,NulLA,NullstellensatzGrobnerChordal,DeLoeraLeeMalkinMargulies_computingInfeasibility_Nullstellensatz,LiLowensteinOmar,DeLoeraLeeMarguliesMiller_Matroids_Polynomials,Pak_Robichaux_2025}. Here we show a new surprising connection of Ramsey numbers to Nullstellensatz certificates.

In Theorem \ref{RamseyNullUpperBound} we give a general construction for Nullstellensatz certificates of Ramsey number upper bounds using the polynomial encoding from Theorem \ref{binEdge2}. In other words, assuming the number of vertices of a Ramsey ideal is $n$ and $n\ge R(G_1,\dots,G_k)$, we give upper bounds for the degrees of the Nullstellensatz certificates. Surprisingly, our construction of 
Nullstellensatz certificates for the Ramsey ideals are related to a certain class of games. 

So-called ``Ramsey games" were studied by Beck in \cite{BeckVDWRamsey}, and he went on to introduce the \emph{online Ramsey numbers} $\tilde R(r,s)$ in \cite{BeckOnlineRamseyProbMethod}, though Kurek and Ruci\'nski independently studied them in \cite{KurekRucinskiOnlineRamsey}. 
These numbers are defined in terms of the following \emph{Builder-Painter game} that is played on an infinite set of vertices. Each turn Builder selects one edge $e$, and Painter selects a color in $[k]$ and colors $e$ with this color. Builder wins the game once a monochromatic $G_i$ in color $i$ is constructed for some $i$. The number $\tilde{R}(G_1,\dots,G_k)$ is the minimum number of edges Builder needs to guarantee a victory. These numbers, as well as some closely related variants, have been studied recently in, for example, \cite{Conlon_Online_Ramsey,DyDzZaOnlinePathCycle,ClemenHeathLavrovOnlinePathCycle,heath2024onlineramseynumbersordered,GuoWarnkeRamseyGamesAlterations,SongWangZhangOnlineK13VsPaths,AdamskiMalgorzataBlavzejOnlineRamseyLongVsShortCycles}.

Our work yields certificate degrees given by the \textit{restricted online Ramsey numbers} $\tilde{R}(G_1,\dots,G_k;n)$, first introduced by Conlon et al. in \cite{RestrictedOnlineRamseyNumbers_Conlon_et_al}. These numbers are defined the same way as their unrestricted counterparts, but the game is played instead on a \emph{finite} set of $n$ vertices. We also use the simplified notation $\tilde R(r_1,\dots,r_k;n)$ for $\tilde R(K_{r_1},\dots,K_{r_k};n)$. It is trivial that $\tilde R(G_1,\dots,G_k) \le \binom{R(G_1,\dots,G_k)}{2}$, but 
it was shown in \cite{RestrictedOnlineRamseyUpperBound} that in the case of 2-color classical Ramsey numbers $\tilde R(r,r;n) \le \binom{n}{2}-\Omega(n\log n)$ for $n = R(r,r)$. While this degree bound is linear in the number of variables, $k\binom n 2$, it is an improvement over the upper bounds for the Nullstellensatz given in, for instance, \cite{BrownawellNullstellensatzBound,Lazard}. For other graphs there are precise results: Briggs and Cox proved that $\tilde{R}(rK_2,rK_2;n) \le n -1$, where $rK_2$ is a matching and $n = R(rK_2, rK_2) = 3r-1$, and Dvo\v{r}\'{a}k later proved the bound is tight \cite{BriggsCoxRestrictedOnlineRamseyMatchingsTrees, RestrictedOnlineRamseyMatchings}. 



\begin{theorem}\label{RamseyNullUpperBound}
	If $n\ge R(G_1,\dots,G_k)$, then there is an explicit Nullstellensatz certificate of degree $\tilde{R}(G_1,\dots, G_k;n)-1$ that the statement $R(G_1,\dots,G_k)> n$ is false using the encoding in Theorem \ref{binEdge2}. In particular, in the case of 2-color classical Ramsey numbers, this implies that if $n \ge R(r,s)$, then there exists a Nullstellensatz certificate of degree $\tilde R(r,s;n)-1$ that the statement $R(r,s) > n$ is false. 
\end{theorem}
\renewcommand{\P}{\mathcal P}

So far we have mostly stated results for graph-theory Ramsey numbers, but
it is very important to stress that Theorem \ref{RamseyNullUpperBound} is just one special case of a more general theorem. The proof of Theorem \ref{RamseyNullUpperBound} does not rely on the graph-theoretic properties of Ramsey numbers, and in fact it applies to the whole of Ramsey theory (arithmetic, geometric, etc).  In particular, we can modify the encoding in Theorem \ref{binEdge2} to suit several well-known problems in Ramsey theory (see, for example, \cite{GrahamRothschildSpencer,LandmanRobertson}). We now express these problems using the general framework below. 

\begin{definition}
	Let $k$ be a positive integer, and let $\{S_n\}$ be a sequence of sets indexed by $n$. For each color $c$ in $[k]$, let $\P_n^c$ be a subset of $S_n$. A triple $A:= (\{S_n\},\{\P_n^c\};k)$ is a Ramsey-type problem if the following hold: 
	\begin{enumerate}
		\item $S_n \subseteq S_{n+1}$ for $i \ge 1$, 
		\item $\P^c_n \subseteq \P^c_{n+1}$ for $n \ge 1$, $1\le c \le k$,
		\item There exists an integer $N$ such that for all $i \ge N$ and every $k$-coloring of $S_i$ there is a color $c$ and some element $X \in \P^c_i$ where each element of $X$ is assigned color $c$. 
	\end{enumerate}
	The smallest such $N$ is called the Ramsey-type number for $A$, and is denoted $R(A)$.
\end{definition}
Note that the subscript $n$ is merely an index and not necessarily equal to $|S_n|$. Moreover, the sets $\P_n^c$ need not be equal for different colors $c$. 

We see that in the problem of computing classical Ramsey numbers $R(r,s)$, we have $k =2$ colors, and the index $n$ denotes the number of vertices in the complete graph whose edges are being 2-colored. Then $S_n = E(K_n) = \{ (i,j) : 1 \le i < j \le n\}$, and the families $\P_n^1$ and $\P_n^2$ consist of all the sets of edges of induced subgraphs of $K_n$ containing $r$ and $s$ vertices, respectively. As another example, the problem of computing \emph{Schur numbers} asks for the smallest $n$ such that every $k$-coloring $[n]$ contains a monochromatic solution to the equation $x+y = z$. In this case we have $S_n = [n]$, and for all $c$ we have $\P_n^c = \{ \{x,y,z\} : \{x,y,z\} \subseteq [n], x+y = z\}$. 


As we see in Section \ref{PolynomialSection}, the encoding in Theorem \ref{binEdge2} can be modified to give bounds for many other Ramsey-type numbers, including Schur, Rado, van der Waerden, and Hales-Jewett numbers. We can define numbers analogous to the restricted online Ramsey numbers for Ramsey-type problems in terms of another Builder-Painter game. We define this game for a fixed $n$ as follows: For each turn, Builder selects one object from $S_n$ and Painter assigns it a color in $[k]$. Builder wins once there is a color $c$ and an element $X \in \P_n^c$ where every element of $X$ is assigned color $c$. Define the number $\tilde R_k(\P_n^1,\dots, \P_n^k;S_n)$ to be the smallest number of turns for which Builder is guaranteed a victory. In this notation, the restricted online Ramsey number $\tilde R(r,s;n)$ is equal to $\tilde R_2(\P_n^1,\P_n^2;S_n)$  with $\P_n$ and $S_n$ defined as above for the Ramsey number $R(r,s)$. Theorem \ref{NullDegree_general_case} generalizes Theorems \ref{binEdge2} and \ref{RamseyNullUpperBound}. 
\begin{theorem}\label{NullDegree_general_case}
	Let $A = (\{S_n\},\{\P^c_n\};k)$ be a Ramsey-type problem. Then for each $n$, the Ramsey-type number for $A$ is strictly greater than $n$ if and only if the following system of equations has no solution over $\Ftwobar$. 
	
	\begin{align*} p_{X,c}:=\prod_{s \in X} x_{c,s} &=0 &&{\forall X \in \P_n^c},\  1\le c \le k,  \\
	q_s:=1+\sum_{i=1}^k x_{i,s}&= 0 &&{\forall s \in S_n}, \\
	u_{i,j,s}:=x_{i,s}x_{j,s} & = 0 &&\forall s\in S_n,\ \forall i,j,\ 1\le i<j\le k.
	\end{align*}
	
	If $n\ge R(A)$, then the minimal degree of a Nullstellensatz certificate for this system is at most $\tilde R_k(\P_n^1,\dots,\P_n^k;S_n) - 1$. 
	
	Moreover, the number of solutions to this system is equal to the number of 
    $k$-colorings of $S_n$ such that for every color $c$, each set $X \in \P_n^c$ contains an object that is not assigned color $c$. 
\end{theorem}


For example, in the case of Schur numbers, the number of solutions to this system is exactly the number of $k$-colorings of $[n]$ that do not contain any monochromatic solutions to $x+y = z$. In Section \ref{PolynomialSection} we give some examples of values of $\tilde R(r,s;n)$ and $\tilde R(\P_n^1,\dots,\P_n^k;S_n)$ and discuss the Nullstellensatz certificates for the associated polynomial systems.

\vskip .3cm
Our second contribution is about  Alon's Combinatorial Nullstellensatz method.
This is a popular technique where combinatorial problems are encoded by a single polynomial $f(x_1,\dots,x_n)$, and the combinatorial property of interest is true depending on whether $f$ vanishes or not at certain testing points. The analysis usually requires finding which coefficients are nonzero. This approach has been used with great success in many situations (see, for example,  \cite{AlonCombinatorialNullstellensatz,GuthPolynomialMethod,ListColoringWithRequests,WongZhu_Graph23Choosable,KarolyiNagy_ProofOfZeilbergerBressoud,BrooksViaAlonTarsi,SauermannWigderson_PolynomialsMultiplicity} and the references therein). 

We show that lower bounds for Ramsey numbers can be obtained by showing that a certain \emph{Ramsey polynomial} $f_{r,s,n}$ is not identically zero. Its coefficients are new combinatorial numbers $E_{n,k,r,H}$ which we call \emph{ensemble numbers}. We show that $E_{n,k,r,H}$ equals the number of ways to choose two distinct edges from $k$-tuples of $r$-cliques inside $K_n$ such that, every edge in a subgraph $H$ is chosen an odd number of times and every edge in its complement $\bar{H}$ is chosen an even number of times. We give a detailed example computing a value of $E_{n,k,r,H}$ in Section \ref{PolynomialSection}.  



Theorem \ref{CombinatorialNullstellensatzTheorem} shows that the numbers $E_{n,k,r,H}$ can be used to find lower bounds for the diagonal Ramsey number $R(r,r)$, and it is an analogue of Theorem 7.2 in \cite{AlonCombinatorialNullstellensatz}. 


\begin{theorem}\label{CombinatorialNullstellensatzTheorem}
	If $$\sum_{k \text{ odd}} 2^k \left(\binom{r}{2}^2-\binom{r}{2}\right)^{\binom{r}{2}-k} E_{n,k,r,H} \neq\sum_{k \text{ even}} 2^k \left(\binom{r}{2}-\binom{r}{2}^2\right)^{\binom{r}{2}-k} E_{n,k,r,H}$$ for some $H$, then $R(r,r)>n$. 
\end{theorem}

\section{Ramsey Numbers and Hilbert's Nullstellensatz} \label{PolynomialSection}

In this section we give several encodings of the problem of finding an upper bound for $R(r,s)$ in terms of the feasibility of a system of polynomial equations.  In the simplest version of the encoding, the variables correspond to edges in the graph $K_n$, and the solutions of the system correspond to graph colorings that avoid monochromatic copies of $K_r$ and $K_s$. If the system is infeasible for some $n$, then no such coloring exists, hence $R(r,s) \le n$. 


Many combinatorial problems can be encoded as a system of polynomial equations, including colorings, independence sets, partitions, etc. (see, e.g., \cite{BussPitassi_InductionPrinciple,BayerThesis,DeLoera_GrobnerBasesGraphColorings,deloera_lee_margulies_onn_2009, NullstellensatzGrobnerChordal,HillarWindfeldt,Lovasz_StableSetsAndPolynomials,MarguliesOnnPasechnik_Nullstellensatz_Partition}).
A Nullstellensatz certificate for such a combinatorial polynomial system is therefore a proof that a combinatorial theorem is true. We are interested on bounding the Nullstellensatz degree for our Ramsey systems.

There are known general ``algebraic geometers'' upper bounds for the degree of a Nullstellensatz certificate, so the above procedure terminates, even when these bounds are exponential and sharp in general \cite{Kollar}. 
However, the exponential bounds should not be bad news for combinatorialists. 
First, it has been shown \cite{Lazard} that for ``combinatorial ideals", the bounds are much better, linear in the number of variables. Over finite fields there are degree bounds that are independent of the number of variables \cite{GreenTao_PolynomialsFiniteFields}, and a recent paper \cite{MoshkowitzYu_FiniteFieldNullstellensatz} gives substantial improvements to these bounds. The bounds we give in Theorems \ref{RamseyNullUpperBound} and \ref{NullDegree_general_case} for our systems of equations are better than the above bounds. Moreover, it has been documented that in practice the degrees of Nullstellensatz certificates of NP-hard problems (e.g., non-3-colorability), tend to be small ``in practice" (see, for example, \cite{MarguliesThesis,LiLowensteinOmar,DeLoeraLeeMalkinMargulies_computingInfeasibility_Nullstellensatz} and the references therein), especially when the polynomial encodings are over finite fields. Note also that when we know the degree of the Nullstellensatz certificate, one can compute explicit coefficients of the Nullstellensatz certificate using a linear algebra system derived by equating the monomials of the identity. This has been exploited in practical computation with great success, see \cite{NulLA,DeLoeraLeeMalkinMargulies_computingInfeasibility_Nullstellensatz,LiLowensteinOmar}. In general, however, certificates are difficult to compute, and for many of the Ramsey problems in this paper it is impractical to compute certificates beyond small cases.

We now prove Theorem \ref{binEdge2} of our encoding for Ramsey numbers over $\Ftwobar$ below. 

\begin{proof}[Proof of Theorem \ref{binEdge2}]
	Suppose there is a solution $\mathbf{x}$ to the system over $\Ftwobar$. Then we color each edge of $K_n$ with the color indicated by the solution vector $\mathbf{x}$. More precisely, for each edge $e$ of $K_n$ and each color $i$, the system has a variable $x_{i,e}$. The polynomials $u_{i,j,e}$ guarantee that for a given $e$, at most one variable $x_{i,e}$ is nonzero. From the polynomials $q_e$, we then see that exactly one index $i$ such that $x_{i,e} = 1$, and let $\phi(\textbf{x})$ be the coloring $\chi$ where $\chi(e)$ is this index. Color each edge $e$ of $K_n$ with the color $\chi(e)$. In the equations involving the polynomials $p_H$, for each subgraph $H$ of $K_n$ with $H \cong G_i$, there is at least one edge $e$ in $H$ with $x_{i,e} = 0$. Therefore $\chi(e) \neq i$, so there is no monochromatic copy of $G_i$ in color $i$. 
	
	Conversely, if we have a coloring $\chi$ of the edges of $K_n$ with no monochromatic $G_i$ in color $i$, then let $\psi(\chi)$ be the solution $\textbf{x}$ where $$x_{i,e}=\begin{cases}
	1 & \text{if } \chi(e) = i, \\
	0 & \text{otherwise}. 
	\end{cases}
	$$ One can check easily that $\textbf{x}$ satisfies the system of equations. 
	The maps $\phi$ and $\psi$ are inverses of each other, and so the number of solutions to the system is equal to the number of colorings of $K_n$ with no monochromatic $G_i$ in color $i$.  
	
	For the first chain of ideals, observe that for a fixed $i$, the polynomial $\prod_{e\in E(H)}x_{i,e}$ divides $\prod_{e\in E(H')}x_{i,e}$ if and only if $H$ is a subgraph of $H'$. Since every copy of $K_{r+1}$ in $K_n$ contains a copy of $K_r$ as a subgraph, in the ideal $RI(n,r+1,s)$, every polynomial of the form $\prod_{e\in E(H')}x_{i,e}$ with $H' \cong K_{r+1}$ is divisible by a generator $\prod_{e\in E(H)}x_{i,e}$ of $RI(n,r,s)$ with $H \cong K_r$. The ideals in the chain are equal for $r > n$ since in this case $K_r$ is not a subgraph of $K_n$. The proof for the second chain of ideals is similar. 
\end{proof}




Before we prove Theorem \ref{RamseyNullUpperBound}, we show a special case as a warm-up example. There is a simple certificate of the fact that $R(r,2) \le r$.

\begin{example} \label{nullR2}
	For all $r$, there exists a Nullstellensatz certificate of degree $\binom r2 -1$ of the statement $R(r,2)\le r$. 
\end{example}

\begin{proof}
	Label the edges of $K_r$ from 1 to $n = \binom{r}{2}$. The following identity is a certificate that $R(r,2)\le r$. Polynomials in parentheses are part of the system of equations in Theorem \ref{binEdge2}.
	
	
	\begin{align*} 1 &= (1+x_{1,1}+x_{2,1})+x_{1,1}(1+x_{1,2}+x_{2,2})+x_{1,1}x_{1,2}(1+x_{1,3}+x_{2,3})+\dots\\ &\quad+x_{1,1}x_{1,2}\cdots x_{1,n-1}(1+x_{1,n}+x_{2,n})\\
	&\quad+x_{2,1}+x_{1,1}(x_{2,2})+x_{1,1}x_{1,2}(x_{2,3})+\dots+x_{1,1}x_{1,2}\cdots x_{1,n-1}(x_{2,n})\\&\quad+(x_{1,1}\cdots x_{1,n})
	\end{align*}
\end{proof}

In the proof of Theorem \ref{RamseyNullUpperBound}, we show how to translate a strategy for Builder into a Nullstellensatz certificate. This method can be used to construct a certificate for all (known) upper bounds for $R(G_1,\dots,G_s)$. Notably, better strategies for Builder yield lower degree certificates. In Example \ref{nullR2}, this is not a concern since the order in which Builder selects edges does not matter, and in fact $\tilde R(r,2) = r$. Painter can simply color every edge Builder selects the first color, and Builder wins only when all $\binom r 2$ edges are selected.

The proofs of Theorems \ref{RamseyNullUpperBound} and \ref{NullDegree_general_case} are similar, and in fact Theorem \ref{RamseyNullUpperBound} follows from Theorem \ref{NullDegree_general_case}, but for the sake of concreteness we begin with Theorem \ref{RamseyNullUpperBound}.




\begin{proof}[Proof of Theorem \ref{RamseyNullUpperBound}]

	Number the edges of $K_n$ from 1 to $\binom{n}{2}$. A $t$-turn game state $g$ is a set $\{(e_{i_1},c_1),(e_{i_2},c_2),\dots,(e_{i_t},c_t)\}$ of pairs of edges $e_{i_j}\in E$ chosen by Builder and colors $c_j\in [k]$ chosen by Painter. A game is complete if there is some color $c \in [k]$ where Painter has colored a monochromatic $G_c$ in color $c$. Let $d:= \tilde R(G_1,\dots,G_k;n)$. If Builder follows an optimal strategy for choosing edges, then the game lasts at most $d$ turns, that is $t\le d$. 
	
	For a $t$-turn game state $g$, define the monomial $\pi(g)$ to be $$\pi(g):= \prod_{j=1}^t x_{c_j,e_{i_j}}.$$ Similarly, for any monomial  $f = \prod_{j=1}^t x_{c_j,e_{i_j}}$ with distinct $e_{i_j}$, let $\sigma(f)$ denote the game state $\{(e_{i_1},c_1),(e_{i_2},c_2),\dots,(e_{i_t},c_t)\}$. 
	
	We will describe an algorithm to construct a Nullstellensatz certificate of the form 
	\begin{equation}\label{Thm2Cert1}\sum_{i=1}^k \sum_{H \cong G_i} \beta_{H,i}p_{H,i} +\sum_{e\in E} \gamma_e q_e = 1. 
	\end{equation}
	
	Denote the left-hand side of Equation \ref{Thm2Cert1} by $L$. For each $i \in [k]$, initialize $\beta_{H,i}$ to 0 for all $H \cong G_i$. Initialize $\gamma_e$ to 0 for all edges $e$ except $e_1$, and set $\gamma_{e_1} := 1$. Then repeat the following:
	
	\begin{enumerate}
		\item Expand and simplify $L$ so that $L$ is a sum of monomials. If $L = 1$, then we are done.
		\item Otherwise, at least one term in $L$ is a nonconstant monomial $f$. 
		\item If $\sigma(f)$ is a completed game state, then $p_{H,i}$ divides $f$ for some color $i$ and $H \cong G_i$. Then set $$\beta_{H,i} \gets \beta_{H,i}+\frac{f}{p_{H,i}}.$$ This results in $L\gets L+f$, which cancels the original $f$ in the certificate since it is an expression over $\Ftwobar$, which has characteristic 2.
		\item If $\sigma(f)$ is not a completed game state, then let $e$ be an edge that Builder should choose in an optimal strategy from the game state $\sigma(f)$. Set $$\gamma_{e}\gets \gamma_{e} + f .$$ Since $fq_e = f+\sum_{i=1}^k fx_{i,e}$, we obtain $L\gets L+f+\sum_{i=1}^k fx_{i,e}$. This results in the cancellation of $f$ in $L$, but adds $k$ additional terms (one for each of Painter's $k$ choices for coloring $e$) to $L$. Note that if $\sigma(f)$ is a $t$-turn game state, then $\sigma(fx_{i,e})$ is a $(t+1)$-turn game state for all $i$.  
	\end{enumerate}
	By the symmetry of $K_n$, it is arbitrary which edge Builder selects first. Therefore each nonconstant term that appears in $L$ corresponds to a game state where Builder (but not necessarily Painter) has followed an optimal strategy. Since terms that correspond to completed games are canceled out in step 3, this procedure terminates, resulting in a Nullstellensatz certificate. Because Builder follows an optimal strategy, the maximal degree of any term in any $\gamma_e$ is $d-1$, so the degree of the certificate is $d-1$.

\end{proof}

To illustrate the importance of Builder's strategy in this method, observe that one can construct a degree 7 certificate for the statement $R(3,3)\le 6$ using the following strategy: For the first five turns, Builder selects each edge incident to some vertex $v$. No matter how Painter colors these edges, three must be colored the same color. Call these edges $vw_1,vw_2,vw_3$. Then for the next three turns, Builder selects the edges $w_1w_2,w_1w_3$, and $w_2w_3$, and Painter must construct a monochromatic triangle (see Figure \ref{FigureBuilderPainter}). However, if Builder plays poorly and selects, for example, the edges $(1,2),(2,3),(3,4),(4,5),(5,6),(1,6),(1,4),(2,5),$ and $(3,6)$, then, no matter what Painter does, there are \emph{no} monochromatic triangles, and this leads to a higher degree certificate.

\begin{figure}[h!]

\caption{Builder-Painter game for $R(3,3)$. For each turn, the edge selected by Builder is indicated in black, and Painter colors it either red or blue.}
\label{FigureBuilderPainter}
\vspace{1em}
\begin{tikzpicture}
    \node[draw] at (-2.5,0) {Turn 1:};
	\draw[fill=black] (0,-1) circle (2pt); 
	\draw[fill=black] (-.866,0.5) circle (2pt);
	\draw[fill=black] (.866,0.5) circle (2pt);
	\draw[fill=black] (0,1) circle (2pt);
	\draw[fill=black] (.866,-0.5) circle (2pt);
	\draw[fill=black] (-.866,-0.5) circle (2pt); 
    \draw[very thick, black] (0,-1) -- (-.866, -0.5);

    \begin{scope}[shift = {(3,0)}]
    \draw[very thick, red] (0,-1) -- (-.866, -0.5);
    \draw[fill=black] (0,-1) circle (2pt); 
	\draw[fill=black] (-.866,0.5) circle (2pt);
	\draw[fill=black] (.866,0.5) circle (2pt);
	\draw[fill=black] (0,1) circle (2pt);
	\draw[fill=black] (.866,-0.5) circle (2pt);
	\draw[fill=black] (-.866,-0.5) circle (2pt);
    
    \end{scope}

      \begin{scope}[shift = {(8.5,0)}]
          \node[draw] at (-2.5,0) {Turn 2:};
 \draw[very thick, red] (0,-1) -- (-.866, -0.5);
 \draw[very thick, black] (0,-1) -- (-.866,0.5);
    \draw[fill=black] (0,-1) circle (2pt); 
	\draw[fill=black] (-.866,0.5) circle (2pt);
	\draw[fill=black] (.866,0.5) circle (2pt);
	\draw[fill=black] (0,1) circle (2pt);
	\draw[fill=black] (.866,-0.5) circle (2pt);
	\draw[fill=black] (-.866,-0.5) circle (2pt);

    \begin{scope}[shift = {(3,0)}]
     
 \draw[very thick, red] (0,-1) -- (-.866, -0.5);
 \draw[very thick, blue] (0,-1) -- (-.866,0.5);
    \draw[fill=black] (0,-1) circle (2pt); 
	\draw[fill=black] (-.866,0.5) circle (2pt);
	\draw[fill=black] (.866,0.5) circle (2pt);
	\draw[fill=black] (0,1) circle (2pt);
	\draw[fill=black] (.866,-0.5) circle (2pt);
	\draw[fill=black] (-.866,-0.5) circle (2pt);
        
    \end{scope}

    \end{scope}

    \begin{scope}[shift = {(0,-3)}]
 \node[draw] at (-2.5,0) {Turn 3:};
	 \draw[very thick, red] (0,-1) -- (-.866, -0.5);
 \draw[very thick, blue] (0,-1) -- (-.866,0.5);
 \draw[very thick, black] (0,-1) -- (0,1);
    \draw[fill=black] (0,-1) circle (2pt); 
	\draw[fill=black] (-.866,0.5) circle (2pt);
	\draw[fill=black] (.866,0.5) circle (2pt);
	\draw[fill=black] (0,1) circle (2pt);
	\draw[fill=black] (.866,-0.5) circle (2pt);
	\draw[fill=black] (-.866,-0.5) circle (2pt);

    \begin{scope}[shift = {(3,0)}]
   \draw[very thick, red] (0,-1) -- (-.866, -0.5);
 \draw[very thick, blue] (0,-1) -- (-.866,0.5);
 \draw[very thick, red] (0,-1) -- (0,1);
    \draw[fill=black] (0,-1) circle (2pt); 
	\draw[fill=black] (-.866,0.5) circle (2pt);
	\draw[fill=black] (.866,0.5) circle (2pt);
	\draw[fill=black] (0,1) circle (2pt);
	\draw[fill=black] (.866,-0.5) circle (2pt);
	\draw[fill=black] (-.866,-0.5) circle (2pt);
    
    \end{scope}

      \begin{scope}[shift = {(8.5,0)}]
          \node[draw] at (-2.5,0) {Turn 4:};
  \draw[very thick, red] (0,-1) -- (-.866, -0.5);
 \draw[very thick, blue] (0,-1) -- (-.866,0.5);
 \draw[very thick, red] (0,-1) -- (0,1);
 \draw[very thick, black] (0,-1) -- (.866,0.5); 
    \draw[fill=black] (0,-1) circle (2pt); 
	\draw[fill=black] (-.866,0.5) circle (2pt);
	\draw[fill=black] (.866,0.5) circle (2pt);
	\draw[fill=black] (0,1) circle (2pt);
	\draw[fill=black] (.866,-0.5) circle (2pt);
	\draw[fill=black] (-.866,-0.5) circle (2pt);

    \begin{scope}[shift = {(3,0)}]
     
  \draw[very thick, red] (0,-1) -- (-.866, -0.5);
 \draw[very thick, blue] (0,-1) -- (-.866,0.5);
 \draw[very thick, red] (0,-1) -- (0,1);
 \draw[very thick, blue] (0,-1) -- (.866,0.5); 
    \draw[fill=black] (0,-1) circle (2pt); 
	\draw[fill=black] (-.866,0.5) circle (2pt);
	\draw[fill=black] (.866,0.5) circle (2pt);
	\draw[fill=black] (0,1) circle (2pt);
	\draw[fill=black] (.866,-0.5) circle (2pt);
	\draw[fill=black] (-.866,-0.5) circle (2pt);
        
    \end{scope}

    \end{scope}

    \end{scope}

    \begin{scope}[shift = {(0,-6)}]
 \node[draw] at (-2.5,0) {Turn 5:};
	\draw[very thick, red] (0,-1) -- (-.866, -0.5);
 \draw[very thick, blue] (0,-1) -- (-.866,0.5);
 \draw[very thick, red] (0,-1) -- (0,1);
 \draw[very thick, blue] (0,-1) -- (.866,0.5); 
 \draw[very thick, black] (0,-1) -- (.866,-0.5);
    \draw[fill=black] (0,-1) circle (2pt); 
	\draw[fill=black] (-.866,0.5) circle (2pt);
	\draw[fill=black] (.866,0.5) circle (2pt);
	\draw[fill=black] (0,1) circle (2pt);
	\draw[fill=black] (.866,-0.5) circle (2pt);
	\draw[fill=black] (-.866,-0.5) circle (2pt);

    \begin{scope}[shift = {(3,0)}]
  \draw[very thick, red] (0,-1) -- (-.866, -0.5);
 \draw[very thick, blue] (0,-1) -- (-.866,0.5);
 \draw[very thick, red] (0,-1) -- (0,1);
 \draw[very thick, blue] (0,-1) -- (.866,0.5); 
 \draw[very thick, red] (0,-1) -- (.866,-0.5);
    \draw[fill=black] (0,-1) circle (2pt); 
	\draw[fill=black] (-.866,0.5) circle (2pt);
	\draw[fill=black] (.866,0.5) circle (2pt);
	\draw[fill=black] (0,1) circle (2pt);
	\draw[fill=black] (.866,-0.5) circle (2pt);
	\draw[fill=black] (-.866,-0.5) circle (2pt);
    
    \end{scope}

      \begin{scope}[shift = {(8.5,0)}]
          \node[draw] at (-2.5,0) {Turn 6:};
   \draw[very thick, red] (0,-1) -- (-.866, -0.5);
 \draw[very thick, blue] (0,-1) -- (-.866,0.5);
 \draw[very thick, red] (0,-1) -- (0,1);
 \draw[very thick, blue] (0,-1) -- (.866,0.5); 
 \draw[very thick, red] (0,-1) -- (.866,-0.5);
 \draw [very thick, black] (-.866,-0.5) -- (0,1);
    \draw[fill=black] (0,-1) circle (2pt); 
	\draw[fill=black] (-.866,0.5) circle (2pt);
	\draw[fill=black] (.866,0.5) circle (2pt);
	\draw[fill=black] (0,1) circle (2pt);
	\draw[fill=black] (.866,-0.5) circle (2pt);
	\draw[fill=black] (-.866,-0.5) circle (2pt);

    \begin{scope}[shift = {(3,0)}]
     
  \draw[very thick, red] (0,-1) -- (-.866, -0.5);
 \draw[very thick, blue] (0,-1) -- (-.866,0.5);
 \draw[very thick, red] (0,-1) -- (0,1);
 \draw[very thick, blue] (0,-1) -- (.866,0.5); 
 \draw[very thick, red] (0,-1) -- (.866,-0.5);
 \draw [very thick, blue] (-.866,-0.5) -- (0,1);
    \draw[fill=black] (0,-1) circle (2pt); 
	\draw[fill=black] (-.866,0.5) circle (2pt);
	\draw[fill=black] (.866,0.5) circle (2pt);
	\draw[fill=black] (0,1) circle (2pt);
	\draw[fill=black] (.866,-0.5) circle (2pt);
	\draw[fill=black] (-.866,-0.5) circle (2pt);
        
    \end{scope}

    \end{scope}

    \end{scope}

\begin{scope}[shift = {(0,-9)}]
 \node[draw] at (-2.5,0) {Turn 7:};
	 \draw[very thick, red] (0,-1) -- (-.866, -0.5);
 \draw[very thick, blue] (0,-1) -- (-.866,0.5);
 \draw[very thick, red] (0,-1) -- (0,1);
 \draw[very thick, blue] (0,-1) -- (.866,0.5); 
 \draw[very thick, red] (0,-1) -- (.866,-0.5);
 \draw [very thick, blue] (-.866,-0.5) -- (0,1);
 \draw [very thick, black] (.866,-0.5) -- (0,1);
    \draw[fill=black] (0,-1) circle (2pt); 
	\draw[fill=black] (-.866,0.5) circle (2pt);
	\draw[fill=black] (.866,0.5) circle (2pt);
	\draw[fill=black] (0,1) circle (2pt);
	\draw[fill=black] (.866,-0.5) circle (2pt);
	\draw[fill=black] (-.866,-0.5) circle (2pt);

    \begin{scope}[shift = {(3,0)}]
  \draw[very thick, red] (0,-1) -- (-.866, -0.5);
 \draw[very thick, blue] (0,-1) -- (-.866,0.5);
 \draw[very thick, red] (0,-1) -- (0,1);
 \draw[very thick, blue] (0,-1) -- (.866,0.5); 
 \draw[very thick, red] (0,-1) -- (.866,-0.5);
  \draw [very thick, blue] (-.866,-0.5) -- (0,1);
 \draw [very thick, blue] (.866,-0.5) -- (0,1);
    \draw[fill=black] (0,-1) circle (2pt); 
	\draw[fill=black] (-.866,0.5) circle (2pt);
	\draw[fill=black] (.866,0.5) circle (2pt);
	\draw[fill=black] (0,1) circle (2pt);
	\draw[fill=black] (.866,-0.5) circle (2pt);
	\draw[fill=black] (-.866,-0.5) circle (2pt);
    
    \end{scope}

      \begin{scope}[shift = {(8.5,0)}]
          \node[draw] at (-2.5,0) {Turn 8:};
   \draw[very thick, red] (0,-1) -- (-.866, -0.5);
 \draw[very thick, blue] (0,-1) -- (-.866,0.5);
 \draw[very thick, red] (0,-1) -- (0,1);
 \draw[very thick, blue] (0,-1) -- (.866,0.5); 
 \draw[very thick, red] (0,-1) -- (.866,-0.5);
 \draw [very thick, blue] (.866,-0.5) -- (0,1);
  \draw [very thick, blue] (-.866,-0.5) -- (0,1);
 \draw [very thick, black] (.866, -0.5) -- (-.866, -0.5); 
    \draw[fill=black] (0,-1) circle (2pt); 
	\draw[fill=black] (-.866,0.5) circle (2pt);
	\draw[fill=black] (.866,0.5) circle (2pt);
	\draw[fill=black] (0,1) circle (2pt);
	\draw[fill=black] (.866,-0.5) circle (2pt);
	\draw[fill=black] (-.866,-0.5) circle (2pt);

    \begin{scope}[shift = {(3,0)}]
     
 \draw[very thick, red] (0,-1) -- (-.866, -0.5);
 \draw[very thick, blue] (0,-1) -- (-.866,0.5);
 \draw[very thick, red] (0,-1) -- (0,1);
 \draw[very thick, blue] (0,-1) -- (.866,0.5); 
 \draw[very thick, red] (0,-1) -- (.866,-0.5);
 \draw [very thick, blue] (.866,-0.5) -- (0,1);
  \draw [very thick, blue] (-.866,-0.5) -- (0,1);
 \draw [very thick, red] (.866, -0.5) -- (-.866, -0.5); 
    \draw[fill=black] (0,-1) circle (2pt); 
	\draw[fill=black] (-.866,0.5) circle (2pt);
	\draw[fill=black] (.866,0.5) circle (2pt);
	\draw[fill=black] (0,1) circle (2pt);
	\draw[fill=black] (.866,-0.5) circle (2pt);
	\draw[fill=black] (-.866,-0.5) circle (2pt);
        
    \end{scope}

    \end{scope}

    \end{scope}
  
\end{tikzpicture} 
\end{figure}

The proof of Theorem \ref{RamseyNullUpperBound} shows that the polynomials can ``simulate" a tree of Builder-Painter games. However, in general the degrees of certificates can be strictly smaller than the bounds given in Theorems \ref{RamseyNullUpperBound}. For example, by a computational search, it can be shown that $\tilde R(3,3;6) = 8$. However, there exists a Nullstellensatz certificate of degree 5 using the encoding in Theorem \ref{binEdge2}, which is better than the bound given in Theorem \ref{RamseyNullUpperBound}. 
We now give the proof of Theorem \ref{NullDegree_general_case}.

\begin{proof}[Proof of Theorem \ref{NullDegree_general_case}]
	Let $A=(\{S_n\}, \{\P_n^c\};k)$ be a Ramsey-type problem. 
	A $t$-turn game state $g$ after $t$ is a set $\{(s_{i_1},c_1),(s_{i_2},c_2),\dots,(s_{i_t},c_t)\}$ of objects $s \in S_n$ chosen by Builder and colors $c_j\in [k]$ chosen by Painter. A game is complete if there is a color $c \in [k]$ and some element $X \in \P_n^c$ where Painter has colored all the elements of $X$ color $c$. Let $d:= \tilde R_k(\P_n;S_n)$. If Builder follows an optimal strategy for choosing edges, then the game lasts at most $d$ turns, that is $t\le d$.
	
	For a game state $g_t$, define the monomial $\pi(g_t)$ to be $$\pi(g_t):= \prod_{j=1}^t x_{c_j,s_{i_j}}.$$ Similarly, for any monomial  $f = \prod_{j=1}^t x_{c_j,s_{i_j}}$ with distinct $s_{i_j}$, let $\sigma(f)$ denote the game state $\{(s_{i_1},c_1),(s_{i_2},c_2),\dots,(s_{i_t},c_t)\}$. 
	
	We will describe an algorithm to construct a Nullstellensatz certificate of the form 
	\begin{equation}\label{Thm2Cert}\sum_{(X,c) \in \P_n} \beta_{X,c}p_{X,c} +\sum_{s\in S_n} \gamma_s q_s = 1. 
	\end{equation}
	
	Denote the left-hand side of Equation \ref{Thm2Cert} by $L$. Initialize $\beta_{X,c}$ to 0 for all $(X,c) \in \P_n$. Let $s^*$ be an object that Builder selects first in an optimal strategy. Initialize $\gamma_{s^*}$ to 1 and $\gamma_s$ to 0 for all other $s \in S_n$. Then repeat the following:
	
	\begin{enumerate}
		\item Expand and simplify $L$ so that $L$ is a sum of monomials. If $L = 1$, then we are done.
		\item Otherwise, at least one term in $L$ is a nonconstant monomial $f$. 
		\item If $\sigma(f)$ is a completed game state, then $p_{X,c}$ divides $f$ for some $(X,c) \in \P_n$. Then set $$\beta_{X,c} \gets \beta_{X,c}+\frac{f}{p_{X,c}}.$$ This results in $L\gets L+f$, which cancels the original $f$ in the certificate since it is an expression over $\Ftwobar$, which has characteristic 2.
		\item If $\sigma(f)$ is not a completed game state, then let $s$ be an object that Builder should choose in an optimal strategy from the game state $\sigma(f)$. Set $$\gamma_{s}\gets \gamma_{s} + f .$$ Since $fq_s = f+\sum_{i=1}^k fx_{i,s}$, we obtain $L\gets L+f+\sum_{i=1}^k fx_{i,e}$. This results in the cancellation of $f$ in $L$, but adds $k$ additional terms (one for each of Painter's $k$ choices for coloring $s$) to $L$. Note that if $\sigma(f)$ is a $t$-turn game state, then $\sigma(fx_{i,s})$ is a $(t+1)$-turn game state for all $i$.  
	\end{enumerate}
	For each nonconstant monomial $f$ that appears in $L$, its corresponding game state $\sigma(L)$ is one where Builder (but not necessarily Painter) has followed an optimal strategy.  Since terms that correspond to completed games are cancelled out in step 3, this procedure terminates and results in a Nullstellensatz certificate. Because Builder follows an optimal strategy, the maximal degree of any term in any $\gamma_s$ is at most $d-1$, so the degree of the certificate is at most $d-1$. 
\end{proof}

As an application of Theorem \ref{NullDegree_general_case}, let $\mathcal{E}$ be a linear equation, and let $R_k(\mathcal E)$ denote the $k$-color \emph{Rado number} for $\mathcal E$, the smallest $n$ such that every $k$-coloring of $[n]$ contains a monochromatic solution to $\mathcal E$. Let $X_{n,\mathcal E}$ be the set of all solutions over $[n]$ to $\mathcal E$. Let $\P_n^c := X_{n,\mathcal E}$ for all $c$. If $R_k(\mathcal E)$ exists, then $(\{[n]\},\{\P_n^c\};k)$ is a Ramsey-type problem, and we have the following corollary. 

\begin{corollary} \label{Corollary_BuilderPainter}
	Let $\mathcal E$ be the linear equation $\sum_{j=1}^t a_i y_i = a_0 $ with a finite Rado number $R_k(\mathcal E)$. Let $X_{n,\mathcal E} = \{(m_1,\dots, m_t): \sum_{j=1}^t a_j m_j = a_0,\ 1\le m_j \le n \}$ be the set of solutions over $[n]$ to $\mathcal E$. Then for every $n$, the following system has no solution over $\Ftwobar$ if and only if $n \ge R_k(\mathcal E)$. 
	
	\begin{align*}
	\prod_{j=1}^t x_{i,m_j} &= 0 & \forall (m_1,\dots, m_t) \in X_{n,\mathcal E},\ 1\le i \le k,\\
	1+ \sum_{i=1}^k x_{i,m} &= 0 & 1\le m \le n, \\
	x_{i,m}x_{j,m} & = 0 & 1 \le m \le n , 1\le i < j \le k. 
	\end{align*}
	
	The degree of a minimal Nullstellensatz certificate for this system has degree at most  $$\tilde R_k(X_{n,\mathcal E},\dots,X_{n,\mathcal E};[n])-1.$$ 
\end{corollary}

As an example, let $\mathcal E$ denote the equation $x+3y = 3z$, and let $X_{9,\mathcal E}$ be the solutions to $\mathcal E$ over $[9]$ as above. It is known that  $R_2(\mathcal E) = 9$ \cite{LandmanRobertson}. However, Builder can select, in order, the integers 4,6,9,3, and 7 to win the Builder-Painter game in at most 5 turns: since $(6,4,6)$ is a solution, 4 and 6 must be different colors, and then since $(9,6,9)$ and $(3,3,4)$ are solutions, 4 and 9 must be one color and 3 and 6 are the other color. But then $(3,6,7)$ and $(9,4,7)$ are solutions, so there is a monochromatic solution no matter which color Painter selects for 7. Therefore $\tilde R_2(X_{9,\mathcal E},X_{9,\mathcal E};[9]) \le 5$, and the minimal degree of a Nullstellensatz certificate for the system of equations in Corollary \ref{Corollary_BuilderPainter} is at most 4. In fact, some computations show the minimal degree is 2. 

Similarly, the encoding in Theorem \ref{binEdge2} for the Schur number $S(2) = R_2(x+y=z)$ also gives an example of Nullstellensatz certificates that are smaller than the ones given by games. It is well-known that $S(2) = 5$, and from the encoding in Theorem \ref{NullDegree_general_case}, we have $S(2) \le 5$ if and only if the following system of equations has no solutions over $\Ftwobar$.
\begin{align*}
&1+x_{1,n}+x_{2,n} = 0, &&1\le n \le 5, \\
&\hspace{18pt} \begin{aligned}
x_{i,1}x_{i,2} &= 0,   &&\hspace{17pt}x_{i,2}x_{i,4} = 0,\\ 
x_{i,1}x_{i,3}x_{i,4}&= 0,   &&x_{i,1}x_{i,4}x_{i,5}= 0, \ \ x_{i,2}x_{i,3}x_{i,5} = 0. \\
\end{aligned}&&\left.\begin{aligned} \\ \\  \end{aligned} \right\rbrace
\text{for} \ i = 1,2.
\end{align*}


A computer search shows that the number $\tilde R_2(X_{5,x+y=z},X_{5,x+y=z};[5]) = 5$, where $X_{5,x+y=z}$ is the set of positive integer solutions to $x+y= z$ in $[1,5]$. The following identity (over $\Ftwobar$) is a degree 3 Nullstellensatz certificate for the above system of equations, which is an improvement on the bound in Theorem \ref{NullDegree_general_case}.  

\begin{align*}
1 &= (x_{2,5}+x_{1,4}x_{1,5}+x_{1,5}x_{2,3}x_{2,4})(1+x_{1,1}+x_{2,1}) + \\
&(x_{1,1}x_{1,3}+x_{2,1}x_{2,5}+x_{1,1}x_{1,5}x_{2,4}+x_{1,1}x_{2,3}x_{2,5}
\\& +x_{1,3}x_{1,5}x_{2,4}+x_{1,4}x_{1,5}x_{2,1})(1+x_{1,2}+x_{2,2}) +\\
&(x_{1,1}x_{2,5}+x_{2,4}x_{1,5}+x_{1,1}x_{1,5}x_{2,4})(1+x_{1,3}+x_{2,3}) +\\
&(x_{1,5}+x_{1,1}x_{1,3}x_{1,5} + x_{1,1}x_{1,3}x_{2,2} +x_{1,2}x_{2,1}x_{2,5})(1+x_{1,4}+x_{2,4})+\\
&(1+x_{1,1}x_{1,3}) (1+x_{1,5}+x_{2,5})+\\
&(x_{1,3}+x_{2,3}x_{2,5}+x_{2,4}x_{1,5})x_{1,1}x_{1,2}+(x_{2,1}x_{1,5}+x_{2,1}x_{2,5}) x_{1,2}x_{2,4} \\
&(x_{2,5}+x_{1,4}x_{1,5})x_{2,1}x_{2,2} +
(x_{1,1}x_{1,3}+x_{1,1}x_{1,5}+x_{1,3}x_{3,5})x_{2,2}x_{2,4} + \\&
(x_{1,5}+x_{2,2})x_{1,1}x_{1,3}x_{1,4}+ \\
&(x_{1,1}x_{1,4}x_{1,5}) +
x_{2,4}(x_{1,2}x_{1,3}x_{1,5})+ 
x_{1,5}(x_{2,1}x_{2,3}x_{2,4})+ \\&
x_{1,2}(x_{2,1}x_{2,4}x_{2,5})+ 
x_{1,1}(x_{2,2}x_{2,3}x_{2,5}).
\end{align*}

{The \emph{van der Waerden number} $w(t,k)$ is the smallest $n$ such that every $k$-coloring of $[n]$ contains a monochromatic $t$-term arithmetic progression \cite{GrahamRothschildSpencer}. Let $AP_{n,t}$ denote the set of all $t$-term arithmetic progressions in $[n]$. Then setting $\P_n^c := AP_{n,t}$ for all $c$, then  $(\{[n]\},\{\P_n^c\};k)$ is a Ramsey-type problem as well.

	\begin{corollary} \label{CorVDW}
		For every $n$, the following system has no solution over $\Ftwobar$ if and only if $n \ge w(t,k)$. 
		\begin{align*}
		\prod_{j=1}^t x_{i,m_j} &= 0 & \forall (m_1,\dots, m_t) \in AP_{n,t},\ 1\le i \le k,\\
		1+ \sum_{i=1}^k x_{i,m} &= 0 & 1\le m \le n, \\
		x_{i,m}x_{j,m} & = 0 & 1 \le m \le n , 1\le i < j \le k. 
		\end{align*}
		
		The minimal degree of a Nullstellensatz certificate for this system is at most $$\tilde R_k(AP_{n,t},\dots,AP_{n,t};[n])-1.$$
		The number of solutions to this system is the number of $k$-colorings of $[n]$ that contain no $t$-term monochromatic arithmetic progressions. 
	\end{corollary}
    
\begin{example}
    It is a well known result that $w(3,2) = 9$, and it can be shown via computer search that $\tilde{R}_2(AP_{9,3},AP_{9,3}) = 6$. Therefore the system of equations in Corollary \ref{CorVDW} for $n = 9, t = 3, k =2$ has a Nullstellensatz certificate of degree at most 5. 
\end{example}
    
	We give one last consequence of Theorem \ref{NullDegree_general_case}. For fixed parameters $t$ and $n$, a \emph{combinatorial line} is a nonconstant sequence of points $v^1,\dots, v^t$, where $v^i \in [t]^n$ such that for every coordinate $j$, the sequence $(v^i_j)_{i=1}^t$ is either constant or $v^i_j = i$ for all $i$. The \emph{Hales-Jewett number $HJ(t,k)$} is the smallest number $n$ such that every $k$-coloring of $[t]^n$ contains a monochromatic combinatorial line \cite{GrahamRothschildSpencer}. Let $L_{t,n}$ denote the set of all combinatorial lines on $t^n$. If we set $\P_n^c = L_{t,n}$ for all $c$, then $(\{[t]^n,\{\P_n^c\};k)$ is a Ramsey-type problem. 
	
	\begin{corollary}
		For every $n$, the following system has no solution over $\Ftwobar$ if and only if $n \ge HJ(t,k)$. 
		\begin{align*}
		\prod_{j=1}^t x_{i,v^j} &= 0 & \forall (v^1,\dots, v^t) \in L_{t,n},\ 1\le i \le k,\\
		1+ \sum_{i=1}^k x_{i,v} &= 0 & v\in [t]^n, \\
		x_{i,v}x_{j,v} & = 0 & v\in[t]^n , 1\le i < j \le k.
		\end{align*}
		The number of solutions to this system is the number of $k$-colorings of $[t]^n$ that do not contain any monochromatic combinatorial lines. 
	\end{corollary}

	\section{Ramsey Numbers and Alon's Combinatorial Nullstellensatz}\label{section_CN}
	
	In this section we introduce a way to encode the problem of finding a lower bound for $R(r,s)$ in terms of properties of a single polynomial. We also define a family of numbers $E_{n,k,r,H}$ whose values can provide bounds for $R(r,r)$. The following theorem of Alon, the ``Combinatorial Nullstellensatz," has been used to solve many problems in combinatorics and graph theory (see, for example, \cite{AlonCombinatorialNullstellensatz,ListColoringWithRequests,BrooksViaAlonTarsi,GuthPolynomialMethod,WongZhu_Graph23Choosable,KarolyiNagy_ProofOfZeilbergerBressoud} and the references therein). 
	\begin{theorem}[Alon, \cite{AlonCombinatorialNullstellensatz}] \label{CN}
		Let $F$ be a field, and let $f\in F[x_1,\dots,x_n]$. Let $\deg(f) = \sum_{i=1}^n t_i$ with each $t_i$ a nonnegative integer, and suppose the coefficient of $\prod_{i=1}^n x_i^{t_i}$ is nonzero. Then if $S_1,\dots, S_n$ are subsets of $F$ with $|S_i|>t_i$, then there exist $s_1\in S_1, \dots, s_n \in S_n$ such that $f(s_1,\dots, s_n) \neq 0$.
	\end{theorem}
	Here we apply the Combinatorial Nullstellensatz to show that lower bounds for Ramsey numbers can be obtained by showing that a certain polynomial is not identically zero. Consider the following polynomial $f(x) = f_{r,s,n}(x)$, where $K_n = (V,E)$ is the complete graph on $n$ vertices:
	
	\begin{equation}\label{CNPoly}
	f(x) =  \left(\prod_{S\subset V,\ |S| = r} \Biggl(\sum_{e\in E(S)} x_e -\binom{r}{2}\Biggr)\right) \left(\prod_{S\subset V,\ |S| = s} \Biggl(\sum_{e\in E(S)} x_e +\binom{s}{2}\Biggr)\right).
	\end{equation}
	
	Every 2-coloring of the edges of $G$ corresponds to an assignment $c:\{x_e\}_{e\in E} \to \{-1,1\}$. If an edge $e$ is colored with the first color, then we set $c(x_e) =1$, and if $e$ is colored with the second color, then $c(x_e) = -1$. Then $f(c(x)) = 0$ if and only if $G$ contains an $r$-clique in the first color or an $s$-clique in the second color. Therefore if $f(c(x_1),\dots,c(x_{|E|})) = 0$ for all colorings $c$, then $R(r,s)\le n$. 
	
	Since we only consider the values of $f$ on $\{-1,1\}^{|E|},$ we may instead consider the multilinear representative of $f$ in the ideal $\left< x_e^2-1\right>_{e\in E}$. This representative can be obtained by deleting each variable with an even exponent from each term in $f$. By Theorem \ref{CN}, this representative is the zero polynomial if and only if $R(r,s) \le |V|$.

	In the proof of Theorem \ref{CombinatorialNullstellensatzTheorem}, we focus on the case when $r=s$ and study the polynomial $f_{r,r,n}$. Before proving Theorem \ref{CombinatorialNullstellensatzTheorem}, we give an example of a value of an ensemble number $E_{n,k,r,H}$, which are coefficients of the multilinear representative of $f_{r,r,n}$. We recall the definition of $E_{n,k,r,H}$ below. 
	
	We call a collection of $k$ $r$-cliques a \emph{$(k,r)$-ensemble}. For each clique in the ensemble, we select exactly two edges, and if each edge in $H$ is selected an odd number of times and each edge of $\bar H$ is selected an odd number of times, then we call this a \emph{valid covering} of a subgraph $H$. The ensemble number $E_{n,k,r,H}$ is the total number of valid coverings of $H$ counted from all $(k,r)$-ensembles of $K_n$.
	
	If $H$ is the graph on five vertices with edge set $\{(1,2),(1,5)\}$, then $E_{5,3,3,H} = 8.$ Figure \ref{FigureEnsembleNumbers} depicts all eight ways.  Each graph has an associated $(3,3)$-ensemble $\mathcal{E}$, and each 3-clique in $\mathcal{E}$ is assigned a distinct color $c\in \{\text{red,green,blue}\}$. The two edges from that 3-clique are colored $c$. Edges colored with more than one color are drawn as multiple edges. In the upper left figure, for example, the associated ensemble is $\{\{1,2,3\},\{1,3,4\},\{1,4,5\} \}$. The edges $(1,2)$ and $(1,3)$ were chosen from the clique $\{1,2,3\}$, the edges $(1,3)$ and $(1,4)$ were chosen from the clique $\{1,3,4\}$, and the edges $(1,4)$ and $(1,5)$ were chosen from the clique $\{1,4,5\}$. The edges of $H$ are chosen exactly once, and all other edges are chosen zero or two times. Note that for some $(3,3)$-ensembles, such as $\{\{1,2,3\},\{1,2,4\},\{3,4,5\}\}$, it is impossible to choose edges in $H$ an odd number of times. 
	

	
	


	

	
	\usetikzlibrary{shapes.geometric}

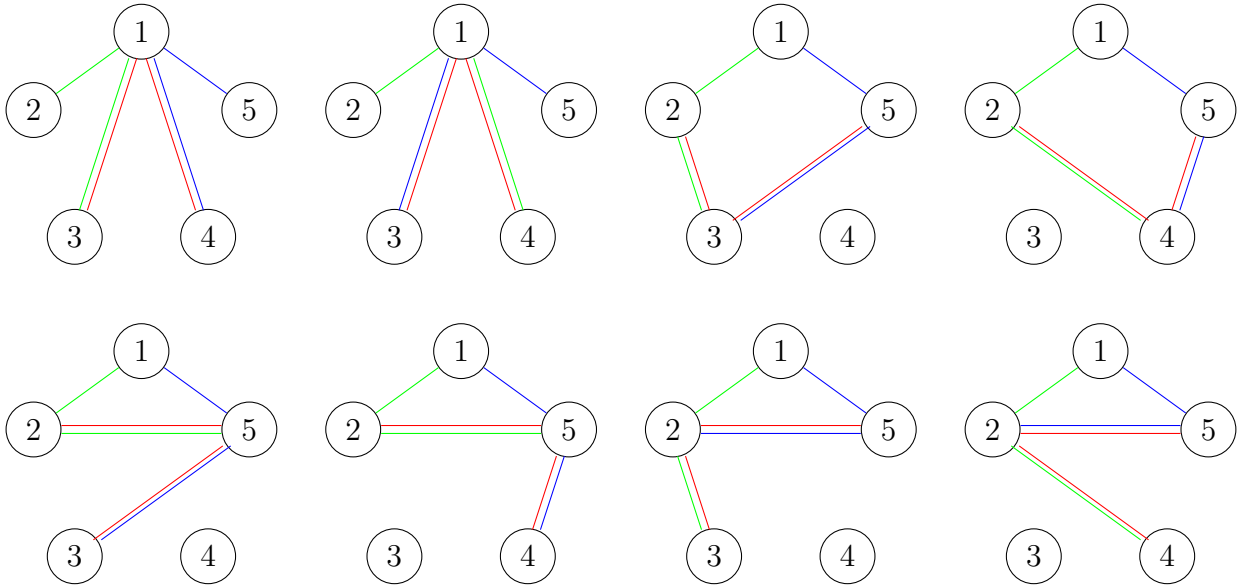
    \begin{figure}
    \caption{Valid coverings of the graph $H$ with edges $(1,2)$ and $(1,5)$.}
    \label{FigureEnsembleNumbers}
	\begin{center}
		\begin{tikzpicture}[scale=0.8, transform shape,mystyle/.style={draw,shape=circle,fill=none},scale = 0.85]
		\def\ngon{5}
		\node[regular polygon,regular polygon sides=\ngon,minimum size=3cm] (p) {};
		\foreach\x in {1,...,\ngon}{\node[mystyle] (p\x) at (p.corner \x){\x};}
		\foreach\x in {1,...,\numexpr\ngon-1\relax}{
			\foreach\y in {\x,...,\ngon}{
			}
		}
		\draw [blue,transform canvas={xshift=1.5pt}] (p1) -- (p4); 
		\draw [red, transform canvas={xshift=-1.5pt}] (p1) -- (p4);
		\draw [blue] (p1) -- (p5); 
		\draw [red, transform canvas={xshift=+1.5pt}] (p1) -- (p3);
		\draw [green] (p1) -- (p2); 
		\draw [green, transform canvas={xshift=-1.5pt}] (p1) -- (p3);
		
		\begin{scope}[shift ={(5,0)}]
		\def\ngon{5}
		\node[regular polygon,regular polygon sides=\ngon,minimum size=3cm] (p) {};
		\foreach\x in {1,...,\ngon}{\node[mystyle] (p\x) at (p.corner \x){\x};}
		\foreach\x in {1,...,\numexpr\ngon-1\relax}{
			\foreach\y in {\x,...,\ngon}{
			}
		}
		\draw [blue,transform canvas={xshift=-1.5pt}] (p1) -- (p3); 
		\draw [red, transform canvas={xshift=-1.5pt}] (p1) -- (p4);
		\draw [blue] (p1) -- (p5); 
		\draw [red, transform canvas={xshift=+1.5pt}] (p1) -- (p3);
		\draw [green] (p1) -- (p2); 
		\draw [green, transform canvas={xshift=1.5pt}] (p1) -- (p4);
		\end{scope}
		\begin{scope}[shift ={(10,0)}]
		\def\ngon{5}
		\node[regular polygon,regular polygon sides=\ngon,minimum size=3cm] (p) {};
		\foreach\x in {1,...,\ngon}{\node[mystyle] (p\x) at (p.corner \x){\x};}
		\foreach\x in {1,...,\numexpr\ngon-1\relax}{
			\foreach\y in {\x,...,\ngon}{
			}
		}
		\draw [blue,transform canvas={xshift=1.5pt}] (p5) -- (p3); 
		
		\draw [blue] (p1) -- (p5); 
		\draw [red, transform canvas={xshift=-1.5pt}] (p5) -- (p3);
		\draw [red, transform canvas={xshift=1.5pt}] (p2) -- (p3);
		\draw [green] (p1) -- (p2); 
		\draw [green, transform canvas={xshift=-1.5pt}] (p2) -- (p3);
		\end{scope}

		\begin{scope}[shift ={(15,0)}]
		\def\ngon{5}
		\node[regular polygon,regular polygon sides=\ngon,minimum size=3cm] (p) {};
		\foreach\x in {1,...,\ngon}{\node[mystyle] (p\x) at (p.corner \x){\x};}
		\foreach\x in {1,...,\numexpr\ngon-1\relax}{
			\foreach\y in {\x,...,\ngon}{
			}
		}
		\draw [blue,transform canvas={xshift=1.5pt}] (p5) -- (p4); 
		
		\draw [blue] (p1) -- (p5); 
		\draw [red, transform canvas={xshift=-1.5pt}] (p5) -- (p4);
		\draw [red, transform canvas={xshift=1.5pt}] (p2) -- (p4);
		\draw [green] (p1) -- (p2); 
		\draw [green, transform canvas={xshift=-1.5pt}] (p2) -- (p4);
		
		\end{scope}
		
		\begin{scope}[shift ={(0,-5)}]
		\def\ngon{5}
		\node[regular polygon,regular polygon sides=\ngon,minimum size=3cm] (p) {};
		\foreach\x in {1,...,\ngon}{\node[mystyle] (p\x) at (p.corner \x){\x};}
		\foreach\x in {1,...,\numexpr\ngon-1\relax}{
			\foreach\y in {\x,...,\ngon}{
			}
		}
		\draw [blue,transform canvas={xshift=1.5pt}] (p5) -- (p3); 
		
		\draw [blue] (p1) -- (p5); 
		\draw [red, transform canvas={xshift=-1.5pt}] (p5) -- (p3);
		\draw [red, transform canvas={yshift=1.5pt}] (p2) -- (p5);
		\draw [green] (p1) -- (p2); 
		\draw [green, transform canvas={yshift=-1.5pt}] (p2) -- (p5);
		
		\end{scope}
		
		\begin{scope}[shift ={(5,-5)}]
		\def\ngon{5}
		\node[regular polygon,regular polygon sides=\ngon,minimum size=3cm] (p) {};
		\foreach\x in {1,...,\ngon}{\node[mystyle] (p\x) at (p.corner \x){\x};}
		\foreach\x in {1,...,\numexpr\ngon-1\relax}{
			\foreach\y in {\x,...,\ngon}{
			}
		}
		\draw [blue,transform canvas={xshift=1.5pt}] (p5) -- (p4); 
		
		\draw [blue] (p1) -- (p5); 
		\draw [red, transform canvas={xshift=-1.5pt}] (p5) -- (p4);
		\draw [red, transform canvas={yshift=1.5pt}] (p2) -- (p5);
		\draw [green] (p1) -- (p2); 
		\draw [green, transform canvas={yshift=-1.5pt}] (p2) -- (p5);
		\end{scope}
		\begin{scope}[shift ={(10,-5)}]

		\def\ngon{5}
		\node[regular polygon,regular polygon sides=\ngon,minimum size=3cm] (p) {};
		\foreach\x in {1,...,\ngon}{\node[mystyle] (p\x) at (p.corner \x){\x};}
		\foreach\x in {1,...,\numexpr\ngon-1\relax}{
			\foreach\y in {\x,...,\ngon}{
			}
		}
		\draw [blue,transform canvas={yshift=-1.5pt}] (p5) -- (p2); 
		
		\draw [blue] (p1) -- (p5); 
		\draw [red, transform canvas={xshift=1.5pt}] (p2) -- (p3);
		\draw [red, transform canvas={yshift=1.5pt}] (p2) -- (p5);
		\draw [green] (p1) -- (p2); 
		\draw [green, transform canvas={xshift=-1.5pt}] (p2) -- (p3);
		\end{scope}
		\begin{scope}[shift ={(15,-5)}]
		\def\ngon{5}
		\node[regular polygon,regular polygon sides=\ngon,minimum size=3cm] (p) {};
		\foreach\x in {1,...,\ngon}{\node[mystyle] (p\x) at (p.corner \x){\x};}
		\foreach\x in {1,...,\numexpr\ngon-1\relax}{
			\foreach\y in {\x,...,\ngon}{
			}
		}
		\draw [blue,transform canvas={yshift=1.5pt}] (p5) -- (p2); 
		
		\draw [blue] (p1) -- (p5); 
		\draw [red, transform canvas={xshift=1.5pt}] (p2) -- (p4);
		\draw [red, transform canvas={yshift=-1.5pt}] (p2) -- (p5);
		\draw [green] (p1) -- (p2); 
		\draw [green, transform canvas={xshift=-1.5pt}] (p2) -- (p4);
		\end{scope}
		
		\end{tikzpicture}
	\end{center}
    \end{figure}
    
	We now prove Theorem \ref{CombinatorialNullstellensatzTheorem}.
	\begin{proof}[Proof of Theorem \ref{CombinatorialNullstellensatzTheorem}]
		We will use the symbol $\equiv$ to denote equivalent representatives in the ideal $I:=\left<x_e^2-1\right>_{e\in E}$. Consider the product of the two terms in $f$ that arise from a fixed $r$-subset $S$ of $V$. Expanding this product and using the relations $x_e^2 \equiv 1$ gives $$\left(\sum_{e\in E(S)}x_e  - \binom r2\right)\left(\sum_{e\in E(S)}x_e  + \binom r2\right) \equiv \left(\sum_{\{e,e'\}\in \binom{E(S)}{2}}2x_ex_{e'} +\binom
		r2 -\binom r2^2\right). $$
		Taking the product of these terms over all $r$-subsets of $V$ gives 
		
		\begin{equation}\label{CNEq2} f(x) \equiv \prod_{S\subset V,\ |S| = r} \left(\sum_{\{e,e'\}\in \binom{E(S)}{2}}2x_ex_{e'} +\binom
		r2 -\binom r2^2\right).
		\end{equation}
		
		After expanding the product, we may write $f$ as a sum of monomials of the form $\prod_{e\in E} x_e^{b_e}$. Two monomials of this form are equivalent modulo $I$ if and only if the parity of $b_e$ is the same for all edges $e$. If $H$ is a subgraph of $G$, it follows that every monomial that satisfies the condition that $b_e$ is odd if and only if $e \in E(H)$ is equivalent to the squarefree monomial $m_H:=\prod_{e\in E(H)}x_e$. Therefore, $f$ is equivalent to a sum of the form\begin{equation}\label{CNEq3}
		f \equiv \sum_{H \subseteq G} a_H m_H.
		\end{equation}
		We now calculate the coefficients $a_H$ in terms of $E_{n,k,r,H}$.
		
		In the expansion of the right-hand side of \eqref{CNEq2}, we have the following combinatorial interpretation of the terms. Each term in the product corresponds to an $r$-subset $S$ of $V$. For each $S$, the term represents a choice of picking either a pair of edges (one of the $2x_ex_{e'}$ terms in the sum), or zero edges (the term $\binom{r}{2} -\binom{r}{2}^2$) from $E(S)$.
		Therefore the coefficient $a_H$ is 
		\begin{equation}\label{CNProof1}a_H = \sum_{k=0}^{\binom n r} 2^k \left(\binom{r}{2}^2-\binom{r}{2}\right)^{\binom{r}{2}-k} E_{n,k,r,H}.\end{equation}
		If $a_H$ is nonzero for some graph $H$, then the multilinear representative of $f$ is nonzero, and by Theorem \ref{CN} there exists a coloring of $K_n$ that makes $f$ nonzero, and in this case it follows that $R(r,r)>n$. Setting the expression $\eqref{CNProof1}$ to be not equal to zero and rearranging terms by the parity of $k$ concludes the proof.  
	\end{proof}
	
	As an example, we give the values of $E_{n,k,r,H}$ for $n=5,r=3$, and $H$ a graph of order five with an even number of edges. We denote by $G_1$ the graph with edge set $\{\{1,2\},\{1,3\},\{2,3\},\{1,4\}\}$ and $G_2$ the graph with edge set $\{\{1,2\},\{2,3\},\{3,4\},\{3,5\}\}$.
	\begin{center}
		\begin{tabular}{c|ccccccccccc}
			\backslashbox{$H$}{$k$} & 0 &1 &2&3&4&5&6&7&8&9&10 \\\hline 
			$\overline {K_5}$ & 1 & 0 & 0 &20&30&132&220&540&585&460&60\\
			$P_3 \cup \overline{K_2}$ &0&1&2&8&44&106&280&496&612&413&86\\
			
			$K_2\cup K_2\cup K_1$ &0&0&4&12&28&124&276&484&628&404&88\\
			
			$K_{1,4}$&0&0&3&4&36&132&242&588&516&428&99 \\
			
			$G_1$& 0&0&2&8&32&120&292&504&592&392&106 \\
			
			$G_2$&0&0&1&10&34&114&292&510&590&390&107 \\
			
			$P_5$&0&0&1&8&40&112&282&520&592&384&109 \\ 
			
			$C_4 \cup K_1$&0&0&2&8&28&136&272&504&612&376&110 \\
			
			$K_2 \cup K_3$&0&0&0&12&36&108&292&516&588&388&108 \\ 
			$\overline{K_{1,4}}$&0&0&0&8&24&120&328&504&552&392&120\\
			
			$\overline{G_1}$& 0&0&0&6&30&118&318&514&554&386&122\\
			
			$\overline{G_2}$&0&0&0&4&36&116&308&524&556&380&124\\
			
			$\overline{P_5}$&  0&0&0&4&32&132&288&524&576&364&128\\
			
			$\overline{C_4\cup K_1}$& 0&0&0&4&36&116&308&524&556&380&124\\
			
			$\overline{K_2 \cup K_3}$& 0&0&0&4&36&108&348&444&636&340&132\\
			
			$\overline{P_3\cup \overline{K_2}}$&0&0&0&0&24&128&344&512&520&384&136\\ 
			
			$\overline{K_2 \cup K_2\cup K_1}$&0&0&0&0&20&144&324&512&540&368&140\\
			
			$K_5$& 0&0&0&0&0&144&400&480&480&400&144 \\
		\end{tabular}
	\end{center}
	Unfortunately, we were unable to find any nontrivial patterns in the above data, and it appears to be difficult to compute the numbers $E_{n,k,r,H}$ in general.

  \subsection*{Acknowledgements} The authors are grateful to Susan Margulies for help with computations and to the anonymous referees for the comments we received. This work was partially supported by NSF grants DMS-1818969, DMS-2348578, and DMS-2434665.

\label{lastpage}

\end{document}